\documentclass[12pt]{article}
\usepackage{enumitem}
\usepackage{graphicx}
\usepackage{tabls}
\usepackage{amsmath}
\usepackage{verbatim}
\usepackage{epsfig}
\usepackage{latexsym}
\usepackage{epic}
\usepackage{eepic}
\usepackage{amssymb}
\usepackage{amsthm}
\usepackage{graphicx}
\usepackage{moreverb}
\usepackage[mathscr]{eucal}
\usepackage[sans]{dsfont}
\usepackage[hidelinks]{hyperref}
\usepackage{url}
\usepackage[margin=1.2in]{geometry} 

\makeatletter
\renewcommand\l@subsection{\@dottedtocline{2}{1.5em}{3em}}
\renewcommand\l@subsubsection{\@dottedtocline{3}{4.5em}{4.5em}}

\newtheorem{Definition}{Definition}[section]
\newtheorem{Theorem}{Theorem}[section]
\newtheorem{Example}{Example}[section]
\newtheorem{Lem}{Lemma}[section]
\newtheorem{Corollary}{Corollary}[section]
\newtheorem{Not}{Note}[section]
\newtheorem{Exer}{Exercise}[section]
\newtheorem{Proposition}{Proposition}[section]

\newcommand{\Ex}{\Example \rm}
\newcommand{\Thm}{\Theorem \rm}

\newcommand{\Cor}{\Corollary \rm}

\renewcommand\footnotesize{%
        \@setfontsize\footnotesize\@ixpt{11}%
        \abovedisplayskip 8\p@ \@plus2\p@ \@minus4\p@
        \abovedisplayshortskip \z@ \@plus\p@
        \belowdisplayshortskip 4\p@ \@plus2\p@ \@minus2\p@
        \def\@listi{\leftmargin\leftmargini
                \topsep 4\p@ \@plus2\p@ \@minus2\p@
                \parsep 2\p@ \@plus\p@ \@minus\p@
                \itemsep \parsep}%
        \belowdisplayskip \abovedisplayskip
}
\makeatother

\begin{document}
\title{\bf Valid confidence intervals for $\mu , \sigma $ 
when there is only one observation available}
\author{Steve Portnoy
\thanks{Corresponding author.}
        \textsubscript{\,}
\thanks{\footnotesize Department of Statistics, University of Illinois at Urbana-Champaign, IL, Email:~\url{sportnoy@illinois.edu} }
\thanks{\footnotesize $\qquad$  and Department of Mathematics and Statistics, Portland, OR}
\\
        Anirban DasGupta
        \thanks{Department of Statistics, Purdue University, West Lafayette, IN. Email:~\url{dasgupta@purdue.edu}.}
}
\date{\Large \today}
\title{\bf \Large Valid confidence intervals for $\mu , \sigma $ when there is only one observation available}
\maketitle

{\large{\sc Abstract}}
Portnoy (2019) considered the problem of constructing an optimal confidence interval 
for the mean based on a single observation $\, X \sim {\cal{N}}(\mu , \, \sigma^2) \,$. Here we extend
this result to obtaining 1-sample confidence intervals for $\, \sigma \,$ and to cases of symmetric
unimodal distributions and of distributions
with compact support. Finally, we extend the multivariate result in Portnoy (2019) to allow a
sample of size $\, m \,$  from a multivariate normal distribution
where $m$ may be less than the dimension.

\bigskip

\newpage
\section{Introduction}
As noted in Portnoy (2019), the problem of constructing a confidence interval 
for the mean based on a single observation $\, X \sim {\cal{N}}(\mu , \, \sigma^2) \,$ 
has been considered for some time. The first published version appears
to be Abbott and Rosenblatt (1962) who showed that the interval 
$\, ( X - c_\alpha | X | , \, X + c_\alpha | X | ) \,$ has coverage probability 
at least $(1 - \alpha)$ for appropriately chosen $c_\alpha$. At about the same
time Charles Stein presented a cloassroom example of the form 
$\, ( -d_\alpha | X | , \, d_\alpha | X | ) \,$, which appeared to have been developed
earlier. Statements attributing
the idea to Herb Robbins (in Rodriguez (1996) and in a
personal communication from Persi Diaconis) suggest that the
example was known to theoretical statisticians before 1960.

Portnoy (2019) considered the problem in somewhat more depth and found a legitimate
confidence interval that is optimal in the sense of having coverage at
least $(1 - \alpha)$ and  minimizing the maximal expected length. This 
interval is in fact randomized,
and it is strictly better than either interval above. This paper also considered a
single observation from a multivariate normal and showed that the set
$\, \{ \, || \mu || \, \leq \, c_\alpha || X ||^2 \, \} \, $ has coverage
probability at least $(1 - \alpha )$ for appropriately chosen $c_\alpha$, and so is a
legitimate confidence set for the multivariate mean.

These results are extended in various directions here. For a single
${\cal{N}}(\mu , \, \sigma^2)$ observation, we find a legitimate confidence interval
for $\sigma$, and we note that these can be combined via the Bonferroni inequality
to form a confidence set for $(\mu, \, \sigma)$ jointly. We also generalize beyond 
the normal distribution to symmetric unimodal distributions and to distributions
with compact support. Finally, we extend the multivariate result to allow a
sample $\, \{ X_1 , \, \cdots \, , \, X_m \} \,$ from the multivariate normal
where $m$ may be less than the dimension, $p$ (and, hence, the covariance matrix
would not be estimable). While the univariate results are mainly an intriguing
curiosity, the multivariate result may actually be useful in some modern large
data problems.

Finally, it may be noted that several probability paradoxes concern possible 
inference based on a single observation. One infamous example is the ``Monte Hall'' 
problem: a prize is placed behind one of three doors, and the contestant may
choose one of the doors. After making a choice, the MC (Monte Hall) shows the 
contestant one of the other two doors which is empty and offers to let the
contestant switch. Since at least 
one of the two doors must be empty, many (most?) people assume the MC 
offered no new information, and take the probability of winning to be the
same (1/2) whether or not the contestant switches. A straightforward calculation
provides the now well-known result that the probability of winning after switching
is 2/3 (assuming all guessing is random).

Another example is the ``two-envelope'' problem: two
players each receive an envelope: one containing the amount X and
the other 2X . By turns, each player may either keep the amount
received or switch envelopes. The conundrum is that if a player 
assumes the envelopes are equally likely, it is always best to switch,
which seems paradoxical. However, having observed one value, the
problem becomes essentially one of hypothesis testing based on a single
observation, and the conditional probabilities will generally fail to
be equally likely. Thus, the player must condition on the observed
value, making the problem one of standard statistical inference, and
not paradoxical. Portnoy (2020) provides a moderately complete treatment
of this and related hypothesis testing problems. See also Wapner (2012)
for a number of such examples.

\section{The Normal Case}

Let $\, X \sim  N(\mu , \sigma ^2) \,$, with $\, -\infty < \mu < \infty \,$, and $\, \sigma > 0$. 
As noted above, there are two classical confidence intervals for the mean, $\mu$, based on the
single observation, $X$: $\, I_1 = I_1(c) \equiv ( X - c | X | , \, X + c | X | ) \,$ (Abbott-Rosenblatt), and 
$\, \{ I_2 = I_2(c) \equiv( -c  | X | , \, c | X | ) \,$ (Stein). The coverage probabilities may be found as direct
corollaries of Theorem 1.1 of Portnoy (2019). The probabilities depend only on $\, \lambda = \mu/\sigma \,$,
and the minimizing values $\lambda^*$ are direct calculations.

\Cor The coverage probability for  $I_1$ depends only on $\lambda$,
is symmetric about zero in $\mu$ (or in $\lambda$), and is given by
\begin{equation}
P_1(\lambda , \, c ) =  \Phi \bigg(\frac{c}{c+1}\, \lambda  \bigg)
+1 - \Phi \bigg (\frac{c}{c-1}\, \lambda \bigg )  \,  . 
\end{equation}
$P_1(\lambda , \, c)$ is minimized over $\lambda$ at 
\begin{equation}
 \lambda^*(c) = \frac{c^2-1}{\sqrt{2}\,c^{3/2}}\,
\sqrt{\log \frac{c+1}{c-1}}. 
\end{equation}
Thus, given $\alpha , 0 < \alpha < 0.5$, there exists a constant $c = c(\alpha )$, such that
$$
\inf_{\mu,\sigma} P_{\mu,\sigma}(I_1(c(\alpha)) \geq 1 - \alpha
$$.

\Cor The coverage probability for  $I_2$ depends only on $\lambda$,
is symmetric about zero in $\mu$ (or in $\lambda$), and is given by 
\begin{equation}
P_2(\lambda , \, c ) =  \Phi \bigg( \frac{c-1}{c}\,\lambda \bigg)
+1 - \Phi \bigg( \frac{c+1}{c}\,\lambda \bigg)  \,  . 
\end{equation}
$P_2(\lambda , \, c)$ is minimized over $\lambda$ at 
\begin{equation}
 \lambda^*(c) = \bigg[ \, \frac{c}{2} \, \log \frac{c+1}{c-1} \, \bigg] ^{1/2} \,\, .
\end{equation}
Thus, given $\alpha , 0 < \alpha < 0.5$, there exists a constant $c = c(\alpha )$, such that
$$
\inf_{\mu,\sigma} P_{\mu,\sigma}(I_2(c(\alpha)) \geq 1 - \alpha
$$.

Note that since the minimal  coverage probability depends only of the length of the interval
(by invariance; see Portnoy (2019)), the minimal coverages for $\, I_1(c) \,$ and $\, I_c(c) \,$ are
exactly the same (as functions of $\, c \,$).

\Thm Given $\alpha , 0 < \alpha <  1$, 
there exists a constant $c = c(\alpha )$, such that 
\[ \inf_{\mu , \sigma }\, P_{\mu , \sigma }\bigg (\sigma^2 \leq \frac{X^2}{c^2}\bigg ) = 1-\alpha . \]
Furthermore, the constant $c = \Phi ^{-1}(\frac{1+\alpha }{2})$, 
where $\Phi (.)$ denotes the standard normal CDF.  

This follows since
$$
P_{\mu , \sigma }\,\bigg (\sigma ^2 \leq \frac{X^2}{c^2}\bigg )
= P_{\mu , \sigma }\,\bigg (\frac{X^2}{\sigma ^2} \geq c^2\bigg)
\geq P_{\mu , \sigma }\,\bigg (\frac{(X-\mu)^2}{\sigma ^2} \geq c^2\bigg)
= P_{\mu , \sigma } ( Z^2 > c^2 )
$$
where the inequality follows since a non-central chi-square has monotone likelihood ratio,
and where $Z$ is a standard unit normal.

\medskip

Finally, note that intervals for $\mu$ and $\sigma$ may be combined (via Bonferroni inequalities) 
to provide a rectangular simultaneous confidence set for $\{\mu , \, \sigma \}$. Of course, other 
simultaneous confidence sets can be constructed, and it seems clear that more circular sets will 
have smaller area. The problem of finding an minimax confidence set (in analogy with the
result in Portnoy (2019)) seems extremely difficult, and it will not be pursued here.

\section{More General Symmetric Unimodal Families} 
\Thm Let $X \sim \frac{1}{\sigma }\, f_0(\frac{x-\mu }{\sigma })$,
where $f_0(-z) = f_0(z)$ for all 
real $z$ and $f_0(.)$ 
is strictly decreasing on $[0, \infty )$. 
Then, given $\alpha , 0 < \alpha < 1$,
there exists $c = c(\alpha )$ such that 
\[ \inf_{\mu , \sigma }\, P_{\mu , \sigma }\bigg (
\sigma \leq \frac{|X|}{c}\bigg ) = 1-\alpha . \]
Furthermore, this constant $c$ satisfies
\begin{equation}
\int_0^c\,f_0(z)\,dz = \frac{\alpha }{2}.
\end{equation} 

\bigskip

\Thm Let $X \sim \frac{1}{\sigma }\, f_0(\frac{x-\mu }{\sigma })$,
where $f_0(.)$ is a continuous 
function and $f_0(-z) = f_0(z)$ for all
real $z$. Let $F_0(.)$ denote the CDF corresponding 
to $f_0, F_0(x) = \int_{-\infty }^x\,f_0(z)dz$.
\\
Let $g(\theta |\alpha ) = \alpha \,f_0(\alpha \,\theta ), 
\alpha , \theta > 0$. Assume further that
\\
(a) $g(\theta |\alpha )$ is strictly MLR in $\theta $.
\\
(b) For all $\alpha _1, \alpha _2$ with $\alpha _2 
> \alpha _1, \log \frac{g(\theta |\alpha _2)}
{g(\theta |\alpha _1)}$ is convex in $\theta $.
\\
Then, 
\\
(i) Given $c > 1$, there exists a unique root
$\theta = \theta (c)$ of the equation
\begin{equation}
\frac{f_0(\frac{c}{c+1}\,\theta )}{f_0(\frac{c}{c-1}\,\theta )}
= \frac{c+1}{c-1}.
\end{equation}
(ii) $\theta (c)$ is continuous in $c$.
\\
(iii) Moreover, for every $c > 1$, 
\begin{equation}
\inf_{\mu , \sigma }\, P\bigg (X-c|X| \leq \mu 
\leq X+c|X|\bigg )
= \psi (c) 
= F_0(\frac{c}{c+1}\,\theta (c))+1-F_0(\frac{c}{c-1}\,\theta (c)).
\end{equation}
(iv) $\psi (c)$ is continuous in $c$. 

\Ex{\bf{The  Cauchy Case:}} Suppose $X \sim C(\mu , \sigma )$,  the Cauchy
distribution with location parameter $\mu $ and scale
parameter $\sigma , -\infty < \mu < \infty , \sigma > 0$.
Therefore, $f_0(z) = \frac{1}{\pi \,(1+z^2)}$,
and direct calculation gives that $\theta (c)$ 
of part (i) of Theorem 3.2 is given by 
$\theta (c) = \frac{\sqrt{c^2-1}}{c}, c > 1$.
It follows that $\theta (c) \to 0$ as $c \to 1$
and $\theta (c) \to 1$ as $c \to \infty $,
and, $\frac{c}{c+1}\,\theta (c) \to 0$
and $\frac{c}{c-1}\,\theta (c) \to \infty $ 
if $c \to 1$, while
both $\frac{c}{c+1}\,\theta (c)$ and
$\frac{c}{c-1}\,\theta (c) \to 1$
if $c \to \infty $. Together, these imply
that $\psi (c)$, the infimum coverage 
probability of part (iii), equation (5),
in Theorem 3.2 satisfies 
$\psi (c) \to 0.5$ as $c \to  1$
and $\psi (c) \to 1$ as $c \to \infty $.
Hence, by the continuity 
of $\psi (c)$ (part (iv), Theorem 3.2)), 
given $\alpha $ such that $0.5 < 1-\alpha < 1$,
there is a $c = c(\alpha )$ 
such that $\psi (c) = 1-\alpha $.
Thus, in the Cauchy case, any nominal 
confidence level $1-\alpha > .5$ 
can be exactly attained by 
a confidence interval of the form $X\pm c|X|$.
\section{General Distributions with Compact Support}
\Thm Let $X \sim F$ and suppose that
$P_F(a \leq |X| \leq b) = 1$, where
$0 < a < b < \infty $. Let $\sigma ^2 
= \sigma^2(F) = \mbox{Var}_F(X)$. \\
Let $K^2 = \frac{4}{(\frac{b}{a}+\frac{a}{b})^2}$, 
and $\alpha > 1-K^2$. Then, 
\[ P_F\bigg (\sigma^2 \leq \frac{X^2}{c^2}\bigg ) \geq 1-\alpha , \]
where $c^2 = 1-\frac{\sqrt{1-\alpha }}{K}$.

\section{A confidence set for $\mu$ based on a sample of size $m$ from 
${\cal{N}}_p(\mu, \, \Sigma)$}

\Thm  \label{MNsamp}
 Let $\, \{ X_1 , \, \cdots \, , \, X_m \} \,$ be a sample from
${\cal{N}}_p(\mu, \, \Sigma)$. Then to achieve
\begin{equation} \label{MNset}
\inf_{\mu , \, \Sigma} P \left\{ || \mu || \leq \frac{c \, || X ||}{\sqrt{m}} \right\} 
\geq 1 - \alpha
\end{equation}
it suffices to take $\, c = 3.85 \, \alpha^{-1/(pm)} \,$.

\section{Proofs} 
{\bf{(Theorem 3.1)}}. If $X \sim \frac{1}{\sigma }\,f_0(\frac{x-\mu }{\sigma })$,
then $Y = \frac{X}{\sigma } \sim f_0(y-\theta )$ 
where $\theta = \frac{\mu }{\sigma }$, and hence,
$|Y| \sim f_0(y-\theta )+f_0(y+\theta )$ under
the assumptions made on $f_0(.)$. Therefore,
\[ P_{\mu , \sigma }\,\bigg (\frac{|X|}{\sigma } \leq c \bigg )
= P_{\mu , \sigma }\,(|Y| \leq c) \]
\[ = \int_0^c \, [f_0(y+\theta )+f_0(y-\theta )]\,dy 
= \int_{\theta }^{\theta +c}\,f_0(z)\,dz+\int_{-\theta }^{c-\theta }\,f_0(z)\,dz \]
\begin{equation} = \int_{\theta -c}^{\theta +c}\,f_0(z)\, dz 
\end{equation}  
(since $f_0(-z) = f_0(z)$ for all z)
\[ \leq \int_{-c}^c\, f_0(z)\, dz = 2\int_{0}^c\, f_0(z)\,dz. \]
(since $f_0(.)$ is strictly decreasing on $(0, \infty )$)
\\
Therefore, if $c$ is chosen such that $\int_{0}^c\, f_0(z)\,dz 
= \frac{\alpha }{2}$, then, we have
\[ P_{\mu , \sigma }\,\bigg (\sigma \geq \frac{|X|}{c}\bigg ) \leq \alpha 
\Rightarrow P_{\mu , \sigma }\,\bigg (\sigma \leq \frac{|X|}{c}\bigg ) \geq 
1-\alpha ,\]
and the infimum of $P_{\mu , \sigma }\,\bigg (\sigma \leq \frac{|X|}{c}\bigg )
= 1-\alpha $ by construction of $c$. This proves Theorem 3.1.
\\\\
{(\bf{Theorem 3.2})}. Following exactly the same lines 
as in Theorem 2.1, one has that
\begin{equation} P_{\mu , \sigma , F_0}\,\bigg (\frac{|X|}{\sigma } \leq c \bigg )
= F_0(\frac{c}{c+1}\,\theta )+1-F_0(\frac{c}{c-1}\,\theta ), 
\end{equation}  
where $\theta = \frac{\mu }{\sigma }$. The minimum 
must be at a critical point, which would satisfy
\[ \frac{c}{c+1}\,f_0(\frac{c}{c+1}\,\theta )-\frac{c}{c-1}\,f_0(
\frac{c}{c-1}\,\theta ) = 0 \]
\begin{equation} \Leftrightarrow \frac{f_0(\frac{c}{c+1}\,\theta )}
{f_0(\frac{c}{c-1}\,\theta )} = \frac{c+1}{c-1}
\Leftrightarrow 
\frac{c}{c+1}\, f_0(\frac{c}{c+1}\,\theta )
= \frac{c}{c-1}\,f_0(\frac{c}{c-1}\,\theta ). \end{equation}  
Since $g(\theta |\alpha ) = \alpha \, f(\alpha \,\theta )$ 
is strictly MLR, it follows that (10) has at most 
one root. However, since 
for any $\alpha _1 < \alpha _2, 
\log \frac{g(\theta |\alpha _2)}{g(\theta |\alpha _1)}$ 
is convex, it follows that
$\frac{f_0(\frac{c}{c+1}\,\theta )}{f_0(\frac{c}{c-1}\,\theta )}
\to \infty $ as $\theta \to \infty $, and
hence (10) must have a root. This establishes 
part (i) of Theorem 3.2. 
\\
Continuity of this unique root, $\theta (c)$ 
follows from joint continuity 
of $\frac{c}{c+1}\, f_0(\frac{c}{c+1}\,\theta )
-\frac{c}{c-1}\,f_0(\frac{c}{c-1}\,\theta )$
in $c$ and $\theta $, as $f_0(z)$ 
has been assumed to be 
continuous in $z$.
The continuity of the infimum $\psi (c)$
follows from continuity of $F_0(.)$ 
and continuity of $\theta (c)$.
\\\\
{(\bf{Theorem 4.1})}. Since $0 \leq a \leq |X| \leq b
< \infty $, by the reverse Cauchy-Schwarz inequality,
\begin{equation} 
E(X^4) \leq \frac{(\frac{b}{a}+\frac{a}{b})^2}{4}\,
[E(X^2)]^2 \Rightarrow 
[E(X^2)]^2 \geq K^2\,E(X^4)
, \end{equation} with $K$ defined 
as in the statement of the theorem.
\\
On the other hand, since $E(X^2) = \mu ^2+\sigma ^2$,
for any $c, 0 < c <1$, by the Paley-Zygmund 
inequality,
\begin{equation}
P(X^2 >c^2\,\sigma ^2) \geq P(X^2 > c^2\,(\mu ^2+\sigma ^2))
\geq (1-c^2)^2\, \frac{[E(X^2)]^2}{E(X^4)}
\geq (1-c^2)^2\,K^2.
\end{equation}
Hence, if $\alpha $ is such that 
$1-\alpha <K^2$, then 
\begin{equation}
P(\sigma ^2 < \frac{X^2}{c^2}) \geq 1-\alpha ,
\end{equation} 
with $c^2$ being chosen as $1-\frac{\sqrt{1-\alpha }}{K}$. 
\\\\
{(\bf{Theorem 5.1})}.
Follow the proof of Theorem 3 in Portnoy (2018) almost exactly. Noting that
$\, \sum_{j=1}^m  || X_j ||^2 \,$ is a non-central chi-square, the 
coverage probability (CP) of the set \eqref{MNset} can be written exactly as in
equation 8 of Portnoy (2018):
\begin{equation} \label{CP1}
CP = P \left\{ || m \, \nu ||^2 \leq c^2 \,  \chi^2_{pm+2K} \right\}
\end{equation}
where $\, K \,$ is Poisson with mean $\, || \nu ||^2 / 2 \,$. Note that the 
only difference here is the appearance of $m$. 

The development and calculations in Portnoy (2018) now go through without any
change except that $p$ is replaced by $mp$. This again provides:
\begin{equation}
1 - CP \leq \alpha \, a^{-pm/2}  \leq \alpha  
\end{equation}
where $\, a = 1 / (1 - \exp ( - 2 \pi \, e^{p/4} + 1 ) ) \,$, and where $c$ 
is defined by
\begin{equation} \label{c2def}
 c^2 = 2 e^2 \, \alpha^{-2/p} \, a \,\, .
\end{equation} 
Theorem \ref{MNsamp} follows trivially by dividing both sides of \eqref{CP1} by $m$.

\end{document}